
\baselineskip=14pt
\parskip=10pt

\magnification=\magstephalf

\def\1{{\overline{1}}}
\def\2{{\overline{2}}}
\parindent=0pt
\overfullrule=0in

\def\frac#1#2{{#1 \over #2}}


\bf
\centerline
{
Automated Proofs (or Disproofs) of Linear Recurrences Satisfied by Pisot Sequences
}

\rm
\bigskip
\centerline
{\it By Shalosh B. EKHAD, N.~J.~A. SLOANE, and  Doron ZEILBERGER}

\bigskip
\centerline
{\it Dedicated to Richard K. GUY on his one hundredth birthday}

{\bf Abstract}: Pisot sequences (sequences $a_n$ with initial terms $a_0=x, a_1=y$, and
defined for $n>1$  by $a_n= \lfloor a_{n-1}^2/a_{n-2} + \frac{1}{2} \rfloor$) often satisfy linear recurrences with constant coefficients
that are valid for all $n \geq 0$, but there are also cautionary examples where
there is a linear recurrence that is valid for an initial range of values of $n$ 
but fails to be satisfied beyond that point,  providing further illustrations of
Richard Guy's celebrated ``Strong Law of Small Numbers''.  In this paper we
present a decision algorithm, fully implemented in
an accompanying  Maple program ({\tt Pisot.txt}), that first
searches for a putative linear recurrence and then decides whether or not 
it holds  for all values of $n$.
We also explain why the failures happen (in some cases 
the `fake' linear recurrence may be valid for thousands of terms).
We conclude by defining, and studying, higher-order analogs of Pisot sequences, 
and point out that
similar phenomena occur there, albeit far less frequently.

{\bf 0. Maple Package and Sample Output}

This article is accompanied by a Maple package, {\tt Pisot.txt}, that is available, along with six input and output files, from 

{\tt http://www.math.rutgers.edu/\~{}zeilberg/mamarim/mamarimhtml/pisot.html} \quad .

{\bf 1. Preface}

Richard Guy famously formulated the {\it Strong Law of Small Numbers}, and in two classic
articles [G1, G2] gave many examples of pairs of sequences that are equal for 
a certain number of initial terms,
 but eventually differ. Before him, around 1820 Charles Babbage [Ba]  had already discussed
numerous examples, and recalled how Fermat was misled by the numbers $2^{2^n}+1$,
and Euler [E] was
almost led to believe that the central trinomial coefficients are the product of
consecutive Fibonacci numbers.

But in all these examples, the sequences only agree for a moderate number of terms.
As shown by David G. Cantor [C1, C2] and David Boyd [B1--B5] , the so-called Pisot sequences provide much
more dramatic examples of Richard Guy's Strong Law of Small Numbers. 
Here we find  pairs of distinct sequences which agree for tens
of thousands of terms.
(Even more extreme examples arise from game theory--see for example
entry {\bf  A078608} in [OEIS], where there are two sequences
which agree for all $n$ from 1 to 777451915729367 but differ at 777451915729368.)

{\bf 2. Pisot Sequences}

We first recall the definition (cf. [Pi], [Ca], [B5]).

{\bf Definition}: The Pisot sequence with index $r$, $E_r(x,y)$ ($0\leq r \leq 1$), where
$0 < x < y$ are integers, is defined by the following nonlinear recurrence:
$$
a_0=x \quad , \quad a_1=y \quad ,
$$
and, for $n>1$,
$$
a_n:= \left\lfloor \, \frac{a_{n-1}^2}{a_{n-2}}  + r \, \right\rfloor \quad ,
$$
where, as usual, $\lfloor \, x \, \rfloor$ denotes the largest integer that is $\leq x$.

The most important special cases are: 

$\bullet$ $r=0$, when $E_0(x,y)$ is abbreviated $T(x,y)$ ,

$\bullet$ $r=\frac{1}{2}$, when $E_{\frac{1}{2}}(x,y)$ is written $E(x,y)$ , and

$\bullet$ $r=1$, when $E_1(x,y)$ is abbreviated $S(x,y)$ .

In the present article we will not consider the limiting cases $r=0$ or $r=1$ (that is, $T(x,y)$ and $S(x,y)$),
although analogous arguments, somewhat more subtle, can be applied to them also.

For many choices of initial conditions $x$, $y$, Pisot sequences do satisfy linear recurrences
that hold for all $n$ (and in this article we present an algorithm--fully implemented in Maple--that
rigorously proves it if this is indeed the case), but there are also many examples where
there exists a recurrence that is valid for a long time, only to eventually break down.

For example, Max Alekseyev  [Al] showed that
$E(5,17)$  (sequence {\bf A010914}) satisfies the linear recurrence
$$
a_n=4 a_{n-1}-2 a_{n-2} \quad ,
$$
for all $n \geq 2$.
On the other hand, David Boyd [B5] found that $E(10,219)$ 
(see sequence {\bf A007699})
satisfies the linear recurrence
$$
a_n=22 a_{n-1}-3 a_{n-2}+18 a_{n-3}-11 a_{n-4} \quad ,
$$
for $4 \leq n \leq 1402$, but that this breaks down at $n=1403$.

Also, one of us (SBE) found (see the bottom of the output file \hfill\break
{\tt http://www.math.rutgers.edu/\~{}zeilberg/tokhniot/oPisot2a.txt}) 
that the Pisot sequence\hfill\break
 $E(30,989)$  ({\bf A276396})
satisfies the recurrence
$$
a_n=33a_{n-1}-2a_{n-2}+30a_{n-3}-11a_{n-4} \quad ,
$$
for $4 \leq n \leq 15888$, but that this breaks down at $n=15889$.

The main tool for explaining  why these Pisot sequences sometimes have such
{\it doppelg\"{a}ngers} (sequences generated by linear recurrences
which agree with them for many terms but eventually differ) is the following result:

{\bf Theorem} (Flor [Fl], Boyd {B5])
{\it If $E_r(x,y)$ ($0 \le r \le 1$) satisfies a linear recurrence then
the defining polynomial  $M(t)$ of the linear recurrence is either $(t-1)^2$ or else
has a single root $r_1 > 1$ outside the unit circle with the remaining roots on or inside the unit circle,
the roots on the unit circle being simple roots.}

As we will see from the analysis and examples below, if there is a second root $r_2$ that is {\it just} outside 
the unit circle, the doppelg\"{a}nger defined by the recurrence can agree with
the Pisot sequence for a large number of terms.

How likely is it that a second root $r_2$  exists outside the unit circle?  If the coefficients of 
the quotient $M(t)/(t-r_1)$ were random (which of course they are not),
then studies of the locations of roots
of random polynomials suggest that the roots tend to be concentrated
in a narrow annulus containing the unit circle (see for example [IZ] 
and the earlier references cited there).  If this were true here then we should expect 
doppelg\"{a}ngers to be fairly common.  Both Cantor [Ca] and Boyd [B1-B5]
have carried out systematic studies of various classes of Pisot sequences.
 It would be nice to have more statistics about 
the minimal polynomials $M(t)$ that arise.

{\bf 3. How to Prove that a Proposed Linear Recurrence for a Pisot Sequence Holds for All Values}

For the sake of pedagogy, before discussing the general case, in this section
we will study a specific example. 

The sequence $E(4,7)$, {\bf A010901}, let's call it $\{ a_n \}$,  starts with
$$
4, 7, 12, 21, 37, 65, 114, 200, 351, 616, 1081, 1897, 3329, 5842, 10252, 17991, 31572, 55405, 97229,  \dots \quad
$$
The OEIS entry formerly  contained the conjecture that this satisfies the linear recurrence
$$
a_n = 2a_{n-1} - a_{n-2} + a_{n-3} \quad ,
$$
wth initial conditions
$$
a_0=4 \quad , \quad a_1=7 \quad , \quad a_2=12 \quad,
$$
together with the remark that this is satisfied for $n \leq 50000$. 
To prove that this holds for all $n$ we proceed as follows
(the same method was used by Max Alekseyev [Al] to establish 
the recurrence for $E(5,17)$ mentioned above).
Recall  that by the definition of Pisot sequences
$$
a_n:= \left\lfloor \, \frac{a_{n-1}^2}{a_{n-2}}  +\frac{1}{2} \, \right\rfloor \quad .
$$

Let's {\bf define} the sequence $b_n$ to be the (obviously unique) sequence satisfying the recurrence
$$
b_n = 2b_{n-1} - b_{n-2} + b_{n-3} \quad ,
$$
subject to the initial conditions
$$
b_0=4 \quad , \quad b_1=7 \quad , \quad b_2=12 \quad .
$$

We have to prove that $a_n=b_n$ for all $n \geq 0$. Using the symmetry of the ``$=$'' relation,
we will prove the equivalent statement that $b_n=a_n$. In other words we must show that
$$
b_n \, = \, \left\lfloor \, \frac{b_{n-1}^2}{b_{n-2}}  +\frac{1}{2} \, \right\rfloor \quad .
$$
But, recalling that $N=\lfloor x \rfloor$ is just shorthand for
$$
N \leq x < N+1 \quad,
$$
our task is to prove that
$$
b_n \leq   \frac{b_{n-1}^2}{b_{n-2}}  +\frac{1}{2} < b_n +1 \quad ,
$$
or equivalently,
$$
- \frac{1}{2} \leq   \, \frac{b_{n-1}^2-b_n b_{n-2}}{b_{n-2}} < \frac{1}{2} \quad .
$$
Define the sequence $c_n$ by
$$
c_n:=b_{n-1}^2-b_n b_{n-2} \quad.
$$

From the linear recurrence defining $b_n$, we know that $b_n$ is given explicitly by
$$
b_n=
3.902586801\,\cdot \, { (1.754877667)}^{n}+ \left(  0.04870659984- 0.09364053397\,i \right)  \left(  0.1225611669+ 0.7448617670\,i \right) ^{n}
$$
$$
+ \left( 
 0.04870659949+ 0.09364053445\,i \right)  \left(  0.1225611669- 0.7448617670\,i \right) ^{n} \quad,
$$
where we have used floating-point  numbers for convenience.
(To make this rigorous we could instead use rational interval
arithmetic. We emphasize that we  do not need to solve the characteristic polynomial of the recurrence exactly,
although in this case of course we could, since it is a cubic polynomial.)

It follows that the sequence $c_n$ is given by
$$
c_n=
0.02472469487\,\cdot \, { (0.5698402912)}^{n}+ \left(  0.4876376523+ 1.233168614\,i \right)  \left(  0.2150798545+ 1.307141279\,i \right)^{n}
$$
$$
+ \left(  0.4876376524- 1.233168615\,i \right)  \left(  0.2150798545- 1.307141279\,i \right)^{n} \quad .
$$
Hence, since the absolute value of the largest terms in $c_n$, 
 $0.2150798545 \pm 1.307141279\,i$, is $1.324717958$, we have
$$
|c_n|=O(1.324717958^n) \quad,
$$
and similarly
$$
b_n = \Omega( 1.754877667^{n}) \quad,
$$
where the implied constants can be easily made explicit if desired. It follows that
$$
\left|\frac{c_n}{b_{n-2}}\right| \, = \, O\left( \left(\frac{1.324717958}{1.754877667}\right)^n\right)=  O((0.7548776664)^n) \quad
$$
and now one can easily find an $N_0$ such that $|\frac{c_n}{b_{n-2}}|< \frac{1}{2}$ for $n \geq N_0$, 
and the
computer can check that this is valid for the first $N_0$ values. This completes the proof.

To get the Maple package to carry out this calculation, you would first load the package by typing {\tt read `Pisot.txt`;} and then running the command

{\tt PtoRv(4, 7, 1/2, 12, 60, 50000);} \quad  .

The arguments to {\tt PtoRv} are  the parameters $x$, $y$, $r$ that define
the Pisot sequence, then the maximal order of a recurrence you wish
to search for (here, 12), then the number of terms of the Pisot sequence $E_r(x,y)$
you would like printed (here, 60),
and finally the number of terms the program should check before giving up (here, 50000).
{\tt PtoRv} is the verbose version; {\tt PtoR} is more succinct.

By using this program we were able to prove conjectured recurrences for 21 entries in the OEIS: 
{\bf A010901, A010904, A010906--A010913, A010924, A020698, A020704, A020720}, \dots .

{\bf 4. The General Case}

Suppose we have found a putative sequence $\{b_n\}$ that appears to agree with a Pisot sequence.
Let $b_n$ satisfy a linear recurrence equation of order $k$  
with constant coefficients, say 
$$
b_n = \sum_{i=1}^{k} A_i b_{n-i}  \quad,
$$
for some integer coefficients $A_1, \dots, A_k$ and given values of  $b_0, \dots, b_{k-1}$.

Let $r_1, \dots, r_k$ be the $k$ roots (for the sake of simplicity we assume that they are distinct) of the
{\bf characteristic polynomial}
$$
t^k-\sum_{i=1}^k A_i t^{k-i} =0 \quad,
$$
and let $r_1$ be the largest root in absolute value, which we assume is real and positive.
(This is reasonable, given the theorem in Section 2.)
Label the roots so that $r_1 >|r_2| \geq |r_3| \geq \dots \geq |r_k|$.

It follows  that $b_n$ satisfies a Binet-type formula
$$
b_n = \sum_{i=1}^{k} C_i r_i^{n} \quad,
$$
for some explicit constants, $C_1, \dots , C_k$ that can easily be  found by linear algebra, 
in terms of the initial values $b_0, \dots, b_{k-1}$.
Hence
$$
c_n:=b_{n-1}^2-b_n b_{n-2}=
\left ( \sum_{i=1}^{k} C_i r_i^{n-1} \right )^2 -
\left ( \sum_{i=1}^{k} C_i r_i^{n} \right )\left ( \sum_{i=1}^{k} C_i r_i^{n-2} \right )
$$
$$
=\, \sum_{i=1}^k \sum_{j=1}^{k} C_i C_j (r_i^{n-1} r_j^{n-1} - r_i^{n} r_j^{n-2})=
\sum_{i=1}^k \sum_{j=1}^{k} C_i C_j r_i^{n-2} r_j^{n-2}(r_i r_j - r_i^2)
$$
$$
= \, \sum_{ { {1\leq i,j \leq k } \atop {i=j}} }  C_i C_j r_i^{n-2} r_j^{n-2}(r_i r_j - r_i^2)+
\sum_{ { {1\leq i,j \leq k } \atop {i \neq j}} }  C_i C_j r_i^{n-2} r_j^{n-2}(r_i r_j - r_i^2)
$$
$$
=\, 0+\sum_{1\leq i<j \leq k}  C_i C_j r_i^{n-2} r_j^{n-2}(2r_i r_j - r_i^2-r_j^2)
$$
$$
= \, -\sum_{1\leq i<j \leq k}  C_i C_j r_i^{n-2} r_j^{n-2} (r_i-r_j)^2 \quad.
$$

Hence $|c_n|=O((r_1 |r_2|)^n)$. We also have that $b_n=\Omega(r_1^n)$.
If $|r_2|<1$, then $\frac{c_n}{b_{n-2}}$ goes to zero exponentially fast, and to check
that
$$
b_n \leq   \frac{b_{n-1}^2}{b_{n-2}}  + r  < b_n +1 \quad ,
$$
once again we need to find an $N_0$ such that for $n \geq N_0$
$$
-r \leq   \frac{c_n}{b_{n-2}} < 1-r \quad ,
$$
and check it for the finitely many cases $n<N_0$.

{\bf 5. Why does E(30,989)'s Doppelg\"{a}nger  Hold for so Many Terms?}

We have already mentioned that the Pisot sequence $E(30,989)$
satisfies the recurrence
$$
a_n=33a_{n-1}-2a_{n-2}+30a_{n-3}-11a_{n-4} \quad ,
$$
for $4 \leq n \leq 15888$ but fails for $n=15889$.

If we apply the above analysis to this recurrence, then
we find that $r_2$ is just {\it outside} the unit circle: $|r_2|=1.00003759711047$, and so
$b_n=a_n$ as long as
$$
 0.2751394860\cdot (1.00003759711047)^n < \frac{1}{2} \quad .
$$
Taking logarithms
$$
(0.00003759629325) \cdot n <  0.5973299074 \quad,
$$
this is true for $n \leq 15888$ but fails beyond that point.

{\bf 6. Infinite Families}

There are many infinite families of Pisot sequences that {\bf do} satisfy linear recurrences. 
Already in 1938 Pisot [Pi] showed that if $x=2$ or $x=3$ then $E(x,y)$
satisfies a linear recurrence of low order, and determined the coefficients.
A very large number of  other families with $x$ in the range 4 to 20  can be viewed here:

{\tt http://www.math.rutgers.edu/\~{}zeilberg/tokhniot/oPisot5.txt} \quad .

For the record here are the first few examples with $x=4, 5, 6$. We denote the unique solution of the linear recurrence (of order $m$)
$$
a_n=\sum_{i=1}^{m} A_i a_{n-i} \quad, \quad a_0=d_1, \dots, a_{m-1}=d_m \quad ,
$$
by the pair of lists 
$$
[[d_1, \dots, d_m] \, , \, [A_1, \dots, A_m]] \quad .
$$

For $k \geq 1$  (and sometimes, if it makes sense, for $k=0$), we have:
$$
{\it E} \left( 4,16\,k+1 \right) =[[4,16\,k+1],[4\,k,k]] \quad ,
$$
$$
{\it E} \left( 4,16\,k+2 \right) =[[4,16\,k+2,64\,{k}^{2}+16\,k+1],[1+4\,k,-2\,k,-k]] \quad ,
$$
$$
{\it E} \left( 4,16\,k+5 \right) =[[4,16\,k+5],[2+4\,k,-1-3\,k]] \quad ,
$$
$$
{\it E} \left( 4,16\,k+7 \right) =[[4,16\,k+7,64\,{k}^{2}+56\,k+12] ,[2+4\,k,-1-k,1+2\,k]] \quad ,
$$
$$
{\it E} \left( 4,16\,k+9 \right) =[[4,16\,k+9,64\,{k}^{2}+72\,k+20] ,[2+4\,k,k,1+2\,k]] \quad ,
$$
$$
{\it E} \left( 4,16\,k+10 \right) =[[4,16\,k+10,64\,{k}^{2}+80\,k+25],[3+4\,k,-2-2\,k,2+3\,k]] \quad ,
$$
$$
{\it E} \left( 4,16\,k+11 \right) =[[4,16\,k+11],[2+4\,k,2+3\,k]] \quad ,
$$
$$
{\it E} \left( 4,16\,k+14 \right) =[[4,16\,k+14,64\,{k}^{2}+112\,k+49],[4+4\,k,-2-2\,k,1+k]] \quad ,
$$
$$
{\it E} \left( 4,16\,k+15 \right) =[[4,16\,k+15],[4+4\,k,-1-k]] \quad ,
$$
$$
\dots \dots \dots
$$
$$
{\it E} \left( 5,25\,k+1 \right) =[[5,25\,k+1],[5\,k,k]] \quad ,
$$
$$
{\it E} \left( 5,25\,k+2 \right) =[[5,25\,k+2,125\,{k}^{2}+20\,k+1],[5\,k,2\,k,k]] \quad ,
$$
$$
{\it E} \left( 5,25\,k+3 \right) =[[5,25\,k+3,125\,{k}^{2}+30\,k+2,625\,{k}^{3}+225\,{k}^{2}+29\,k+1],[5\,k,3\,k,2\,k,k]] \quad .
$$
$$
\dots \dots \dots
$$
$$
{\it E} \left(6, 36\,k + 1 \right) = [[6, 36\,k + 1], [6\,k, k]] \quad ,
$$
$$
{\it E} \left(6, 36\,k + 2 \right) = [[6, 36\,k + 2, 216\,k^2 + 24\,k + 1], [6\,k, 2\,k, k]] \quad ,
$$
$$
{\it E} \left(6, 36\,k + 3 \right) = [ [6, 36\,k + 3, 216\,k^2  + 36\,k + 2, 1296\,k^3  + 324\,k^2  + 33\,k + 1], [1 + 6\,k, -3\,k, -k, -k]] \quad .
$$
$$
\dots \dots \dots
$$

Note that Pisot [P1], Cantor [C2], Boyd [B5] already observed that
the Pisot sequences  $E_r(x,y)$ tend to form families whose properties 
depend on the value of $y$ mod $x^2$, That is, the sequences
$E_r(x, kx^2+j)$, $k=0,1,2,\dots$ all tend to satisfy similar linear recurrences, or 
appear not to satisfy such a recurrence.  The above examples are consistent with this
observation.  

It is likely that some of our results for x = 4 and 5 were already known to 
Galyean [Ga], but we have not been able to get access to his dissertation. 

It is also possible to find {\it doubly}-infinite (i.e. two-parameter) families, but we stop here.


{\bf 7. Higher-Order Generalizations}

A crucial property of Pisot sequences  is that $a_n a_{n+2} - a_{n+1}^2$ is small compared to $a_n$.
Since
$$
a_n a_{n+2} - a_{n+1}^2 = \det \, \pmatrix{ a_n & a_{n+1} \cr a_{n+1} & a_{n+2} } \quad,
$$
it is natural to generalize the definition, and to consider sequences for which, for some $s>1$, the Hankel determinant
$$
\Delta_s :=
\det \,  \pmatrix{ a_n & \dots & a_{n+s} \cr a_{n+1} &  \dots  & a_{n+s+1}  \cr  \dots & \dots & \dots \cr  \dots 
\cr \dots & \dots & \dots \cr  a_{n+s} & \dots & a_{n+2s }} 
$$
is small.

Note that for any sequence that satisfies a linear recurrence with constant coefficients of
order $s$,  the above determinant is identically zero.

Let us define $F_s$ and $G_s$ by writing
$$
\Delta_s= a_{n+2s} F_s(a_1, \dots, a_{n+2s-1}) - G_s(a_1, \dots, a_{n+2s-1})  \quad .
$$
Then we  define an order-$s$ Pisot sequence, $E_r(a_0, \dots , a_{2s-1})$ with parameter  $r$ ($0 \leq r \leq 1$)
by the rules that for $0\leq n \leq 2s-1$ the value is $a_n$, and for $n \geq 0$ we have
$$
a_{n+2s}=\left\lfloor \,   \frac{G_s(a_1, \dots, a_{n+2s-1})}{ F_s(a_1, \dots, a_{n+2s-1})}  +r \, \right\rfloor \quad .
$$

A calculation analogous to that in Section 4
shows that a necessary condition for a linear recurrence  with constant coefficients
to be an order-$s$ generalized Pisot sequence is that the $(s+1)$-st largest absolute value of the roots is less then 1. (Presumably there is also an analog of the theorem
in Section 2 which applies here.)
See the output file

{\tt http://www.math.rutgers.edu/~zeilberg/tokhniot/oPisot4.txt} \quad ,

for numerous examples.

{\bf 8. References}

[Al] Max Alekseyev, {\it Comments on Sequence} {\bf A010914}, Sep. 3 2013; {\tt https://oeis.org/A010914}.

[Ba] Charles Babbage, {\it Of induction}, from {\it Essays in the Philosophy of Analysis}, circa 1820,  typeset and edited by
Martin Fagereng Johansen, 2013; available from  \hfill\break
{\tt http://martinfjohansen.com/ofinduction/ofinduction-2013-11-11-firstdraft.pdf} \quad .

[B1] David W. Boyd, {\it Pisot sequences which satisfy no linear recurrences}, Acta Arith., {\bf 32} (1977) 89--98.

[B2] David W. Boyd, {\it Some integer sequences related to the Pisot sequences}, Acta Arith., {\bf 34} (1979), 295--305.

[B3] David W. Boyd, {\it On linear recurrence relations satisfied by Pisot sequences}, Acta Arith., {\bf 47} (1986) 13--27; {\bf 54} (1990), 255--256.

[B4] David W. Boyd, {\it Pisot sequences which satisfy no linear recurrences, II}, Acta Arith. {\bf 48} (1987) 191--195.

[B5] David W. Boyd, {\it Linear recurrence relations for some generalized Pisot sequences},
in {\it Advances in Number Theory (Kingston ON, 1991)}, 
pp. 333--340, Oxford Univ. Press, New York, 1993; with updates from 1996 and 1999;
available from \hfill\break
https://www.researchgate.net/profile/David\_Boyd7/publication/262181133\_Linear\_recurrence\_relations\_for\_some\_generalized\_Pisot\_sequences\_-\_annotated\_with\_corrections\_and\_additions/links/00b7d536d49781037f000000.pdf \quad .

[C1] David G. Cantor, 
{\it Investigation of T-numbers and E-sequences}, in
A.~O.~L.~ Atkins and B.~J.~Birch, eds., 
{\it Computers in Number Theory},   Acad. Press, NY,  1971, pp. 137--140.

[C2] David G. Cantor, {\it On families of Pisot} $E${\it -sequences}, 
Ann. Scient. \'{E}c. Norm. Sup.,
 {\bf 9} (1976),  283--308.

[E] Leonhard Euler, {\it Exemplum memorabile inductionis fallacis}, Opera Omnia, Series Prima, 
{\bf 15} (1911), 50--69, Teubner. Leipzig, Germany.

[Fl] Peter Flor, 
{\it \"{U}ber eine Klasse von Folgen nat\"{u}rlicher Zahlen}, 
Math. Annalen, {\bf 140} (1960), 299--307.

[Ga] Paul H. Galyean,
{\it On Linear Recurrence Relations for E-Sequences},
Ph.D. Dissertation, Univ. California Los Angeles, 1971 (unpublished).

[G1] Richard K. Guy, {\it The strong  law of small numbers}, Amer. Math. Monthly, {\bf 95} (1988), 697--712.

[G2] Richard K. Guy, {\it The second strong  law of small numbers}, Math. Mag. {\bf 63} (1990), 3--20.

[IZ]
Ildar Ibragimov and Dmitry Zaporozhets, 
{\it On distribution of zeros of random polynomials in complex plane}, 
in {\it Prokhorov and Contemporary Probability Theory}, 
Springer, 2013, pp.~303--323; available from arXiv:1102.3517..

[OEIS] 
The OEIS Foundation Inc.,
{\it The On-Line Encyclopedia of Integer Sequences}, 
{\tt https://oeis.org}.

[Pi]
Charles Pisot,
{\it La r\'{e}partition modulo $1$ et les nombres alg\'{e}briques},
Ann. Scuola Norm. Sup. Pisa Cl. Sci.,
{\bf 7} (1938), 205--248.

\bigskip
\bigskip

\hrule
\bigskip
Shalosh B. Ekhad, c/o D. Zeilberger, Department of Mathematics, Rutgers University (New Brunswick), Hill Center-Busch Campus, 110 Frelinghuysen
Rd., Piscataway, NJ 08854-8019, USA.
\bigskip
\hruleÄ
\bigskip
N. J. A. Sloane, The OEIS Foundation Inc, 11 South Adelaide Ave, Highland Park, NJ 08904, USA,  and Department of Mathematics, Rutgers University (New Brunswick);  \hfill \break
 njasloane at gmail dot com \quad ;  \quad {\tt  http://neilsloane.com/ } \quad .
\bigskip
\hrule
\bigskip
Doron Zeilberger, Department of Mathematics, Rutgers University (New Brunswick), Hill Center-Busch Campus, 110 Frelinghuysen
Rd., Piscataway, NJ 08854-8019, USA; \hfill \break
DoronZeil at gmail dot com \quad ;  \quad {\tt http://www.math.rutgers.edu/\~{}zeilberg/} \quad .
\bigskip
\hrule
\bigskip
\bigskip
To be published in The Personal Journal of Shalosh B. Ekhad and Doron Zeilberger  \hfill \break
({\tt http://www.math.rutgers.edu/\~{}zeilberg/pj.html}), 
N.~J.~A.~Sloane's Home Page \hfill \break ({\tt  http://neilsloane.com/}),
and the {\tt arxiv.org} \quad .
\bigskip
\bigskip
\hrule
\bigskip

{\bf Sept. 20, 2016}

\end